\documentclass[]{amsart}
\usepackage{amssymb,amsmath}
\usepackage[mathscr]{euscript}

\newcounter{sec}

\newcounter{punct}[sec]

\def\punct{\refstepcounter{punct}{\arabic{sec}.\arabic{punct}.  }}

\newtheorem{theorem}{Theorem}[sec]
\newtheorem{proposition}[theorem]{Proposition}

\newtheorem{lemma}[theorem]{Lemma}

\newtheorem{corollary}[theorem]{Corollary}
\newtheorem{observation}[theorem]{Observation}

\def\COUNTERS{\addtocounter{sec}{1}
              \setcounter{punct}{0}
          \setcounter{equation}{0}
          \setcounter{theorem}{0}
         }
          
          \def\sm{\smallskip}
          
\begin{document}

\newcommand{\supp}{\mathop {\mathrm {supp}}\nolimits}
\newcommand{\rk}{\mathop {\mathrm {rk}}\nolimits}
\newcommand{\Aut}{\mathop {\mathrm {Aut}}\nolimits}
\newcommand{\Out}{\mathop {\mathrm {Out}}\nolimits}
\newcommand{\OO}{\mathop {\mathrm {O}}\nolimits}
\renewcommand{\Re}{\mathop {\mathrm {Re}}\nolimits}
\newcommand{\ch}{\cosh}
\newcommand{\sh}{\sinh}

\def\0{\mathbf 0}

\def\ov{\overline}
\def\wh{\widehat}
\def\wt{\widetilde}

\renewcommand{\rk}{\mathop {\mathrm {rk}}\nolimits}
\renewcommand{\Aut}{\mathop {\mathrm {Aut}}\nolimits}
\renewcommand{\Re}{\mathop {\mathrm {Re}}\nolimits}
\renewcommand{\Im}{\mathop {\mathrm {Im}}\nolimits}
\newcommand{\sgn}{\mathop {\mathrm {sgn}}\nolimits}
\newcommand{\Isoc}{\mathop {\mathrm {Isoc}}\nolimits}
\newcommand{\PIsoc}{\mathop {\mathrm {PIsoc}}\nolimits}

\newcommand{\Sch}{\mathop {\mathrm {Sch}}\nolimits}
\newcommand{\sch}{\mathop {\mathrm {sch}}\nolimits}
\newcommand{\Fr}{\mathop {\mathrm {Fr}}\nolimits}
\newcommand{\Op}{\mathop {\mathrm {Op}}\nolimits}
\newcommand{\Mat}{\mathop {\mathrm {Mat}}\nolimits}
\newcommand{\Exp}{\mathop {\mathrm {Exp}}\nolimits}

\def\bfa{\mathbf a}
\def\bfb{\mathbf b}
\def\bfc{\mathbf c}
\def\bfd{\mathbf d}
\def\bfe{\mathbf e}
\def\bff{\mathbf f}
\def\bfg{\mathbf g}
\def\bfh{\mathbf h}
\def\bfi{\mathbf i}
\def\bfj{\mathbf j}
\def\bfk{\mathbf k}
\def\bfl{\mathbf l}
\def\bfm{\mathbf m}
\def\bfn{\mathbf n}
\def\bfo{\mathbf o}
\def\bfp{\mathbf p}
\def\bfq{\mathbf q}
\def\bfr{\mathbf r}
\def\bfs{\mathbf s}
\def\bft{\mathbf t}
\def\bfu{\mathbf u}
\def\bfv{\mathbf v}
\def\bfw{\mathbf w}
\def\bfx{\mathbf x}
\def\bfy{\mathbf y}
\def\bfz{\mathbf z}

\def\bfA{\mathbf A}
\def\bfB{\mathbf B}
\def\bfC{\mathbf C}
\def\bfD{\mathbf D}
\def\bfE{\mathbf E}
\def\bfF{\mathbf F}
\def\bfG{\mathbf G}
\def\bfH{\mathbf H}
\def\bfI{\mathbf I}
\def\bfJ{\mathbf J}
\def\bfK{\mathbf K}
\def\bfL{\mathbf L}
\def\bfM{\mathbf M}
\def\bfN{\mathbf N}
\def\bfO{\mathbf O}
\def\bfP{\mathbf P}
\def\bfQ{\mathbf Q}
\def\bfR{\mathbf R}
\def\bfS{\mathbf S}
\def\bfT{\mathbf T}
\def\bfU{\mathbf U}
\def\bfV{\mathbf V}
\def\bfW{\mathbf W}
\def\bfX{\mathbf X}
\def\bfY{\mathbf Y}
\def\bfZ{\mathbf Z}
\def\bfT{\mathbf T}

\def\frF{\mathfrak F}
\def\frD{\mathfrak D}
\def\frX{\mathfrak X}
\def\frS{\mathfrak S}
\def\frZ{\mathfrak Z}
\def\frL{\mathfrak L}
\def\frG{\mathfrak G}
\def\frg{\mathfrak g}
\def\frh{\mathfrak h}
\def\frf{\mathfrak f}
\def\frl{\mathfrak l}
\def\frp{\mathfrak p}
\def\frq{\mathfrak q}
\def\frr{\mathfrak r}

\def\bfw{\mathbf w}

\def\R {{\mathbb R }}
 \def\C {{\mathbb C }}
  \def\Z{{\mathbb Z}}
\def\K{{\mathbb K}}
\def\N{{\mathbb N}}
\def\Q{{\mathbb Q}}
\def\A{{\mathbb A}}
\def\U{{\mathbb U}}

\def\T{\mathbb T}
\def\P{\mathbb P}

\def\G{\mathbb G}

\def\cD{\EuScript D}
\def\cL{\mathscr L}
\def\cK{\EuScript K}
\def\cM{\EuScript M}
\def\cN{\EuScript N}
\def\cP{\EuScript P}
\def\cQ{\EuScript Q}
\def\cR{\EuScript R}
\def\cT{\EuScript T}
\def\cW{\EuScript W}
\def\cY{\EuScript Y}
\def\cF{\EuScript F}
\def\cG{\EuScript G}
\def\cZ{\EuScript Z}
\def\cI{\EuScript I}
\def\cB{\EuScript B}
\def\cA{\EuScript A}
\def\cO{\EuScript O}
\def\cE{\EuScript E}
\def\ex{\EuScript E}

\def\bbA{\mathbb A}
\def\bbB{\mathbb B}
\def\bbD{\mathbb D}
\def\bbE{\mathbb E}
\def\bbF{\mathbb F}
\def\bbG{\mathbb G}
\def\bbI{\mathbb I}
\def\bbJ{\mathbb J}
\def\bbL{\mathbb L}
\def\bbM{\mathbb M}
\def\bbN{\mathbb N}
\def\bbO{\mathbb O}
\def\bbP{\mathbb P}
\def\bbQ{\mathbb Q}
\def\bbS{\mathbb S}
\def\bbT{\mathbb T}
\def\bbU{\mathbb U}
\def\bbV{\mathbb V}
\def\bbW{\mathbb W}
\def\bbX{\mathbb X}
\def\bbY{\mathbb Y}

\def\kappa{\varkappa}
\def\epsilon{\varepsilon}
\def\phi{\varphi}
\def\le{\leqslant}
\def\ge{\geqslant}

\def\B{\mathrm B}

\def\la{\langle}
\def\ra{\rangle}

\def\lambdA{{\boldsymbol{\lambda}}}
\def\alphA{{\boldsymbol{\alpha}}}
\def\betA{{\boldsymbol{\beta}}}
\def\mU{{\boldsymbol{\mu}}}

\def\const{\mathrm{const}}
\def\rem{\mathrm{rem}}
\def\even{\mathrm{even}}
\def\SO{\mathrm{SO}}
\def\SL{\mathrm{SL}}
\def\SU{\mathrm{SU}}
\def\GL{\operatorname{GL}}
\def\End{\operatorname{End}}
\def\Mor{\operatorname{Mor}}
\def\Aut{\operatorname{Aut}}
\def\inv{\operatorname{inv}}
\def\red{\operatorname{red}}
\def\Ind{\operatorname{Ind}}
\def\dom{\operatorname{dom}}
\def\im{\operatorname{im}}
\def\md{\operatorname{mod\,}}
\def\St{\operatorname{St}}
\def\Ob{\operatorname{Ob}}
\def\PB{{\operatorname{PB}}}
\def\Tra{\operatorname{Tra}}

\def\ZZ{\mathbb{Z}_{p^\mu}}
\def\F{\mathbb{F}}

\def\cH{\EuScript{H}}
\def\cQ{\EuScript{Q}}
\def\cL{\EuScript{L}}
\def\cX{\EuScript{X}}

\def\Di{\Diamond}
\def\di{\diamond}

\def\fin{\mathrm{fin}}
\def\ThetA{\boldsymbol {\Theta}}

\def\0{\boldsymbol{0}}

\def\FF{\,{\vphantom{F}}_3F_2}
\def\HH{\,\vphantom{H}^{\phantom{\star}}_3 H_3^\star}
\def\Ho{\,\vphantom{H}_2 H_2}

\def\disc{\mathrm{disc}}
\def\cont{\mathrm{cont}}

\def\fan{\vphantom{|^|}}

\def\osigma{\ov\sigma}
\def\ot{\ov t}

\def\Afr{\mathrm{Afr}}
\def\fr{\mathfrak{fr}}
\def\Fr{\mathrm{Fr}}

\def\tri{|\!|\!|}

\begin{center}
\Large\bf
A Lie group corresponding to the free Lie algebra
\\
 and its universality

\bigskip

\large\sc Yury A. Neretin%
\footnote{The work is supported by the grant of FWF (Austrian Scientific Funds), PAT5335224.}
\end{center}
\vspace{22pt}

{\small
Consider the real free Lie algebra $\mathfrak{fr}_n$ with generators 
$\omega_1$, \dots, 
$\omega_n$. Since it is positively graded, it has a
completion $\overline{\mathfrak{fr}}_n$ consisting of formal series.
By the Campbell--Hausdorff formula, we have a corresponding Lie group 
$\overline{\mathrm{Fr}}_n$. It is the set 
$\exp\bigl(\overline{\mathfrak{fr}}_n\bigr)$ in the completed universal
 enveloping algebra of $\mathfrak{fr}_n$.
Also, the group $\overline{\mathrm{Fr}}_n$ is a  'submanifold' in the algebra of formal
associative noncommutative  series 
in $\omega_1$, \dots, 
$\omega_n$, the 'submanifold' is determined by a certain system of quadratic equations. We consider a certain dense subgroup 
$\mathrm{Fr}_n^\infty\subset \overline{\mathrm{Fr}}_n$
with a stronger (Polish)
 topology and show that any 
 homomorphism $\pi$ from $\mathfrak{fr}_n$ to a real finite-dimensional Lie algebra
 $\mathfrak{g}$ can be integrated in a unique way to a  homomorphism $\Pi$ from $\mathrm{Fr}_n^\infty$ to the corresponding simply connected Lie group $G$.
If $\pi$ is surjective, then $\Pi$ also is surjective. 
  Note that
 Pestov (1993) constructed a separable Banach--Lie group such that
 any separable Banach--Lie group is its quotient.

{\bf Key words.} Free Lie algebra, universal enveloping algebra, Lie groups,  Banach algebra, Banach--Lie group.

}

\section{The group $\ov\Fr_n$. Preliminaries}

\COUNTERS

{\bf \punct Free associative algebras.}
Denote by $\Afr_n$ the real free  associative algebra
generated by $\omega_1, \dots , \omega_n$.

Denote by $\cA$ the set of all words $\alpha=\alpha_1\dots\alpha_k$,
where $k=0$, $1$, $2$, $\dots$, and each $\alpha_j$ is a positive
integer $\le n$.
Denote
$$
|\alpha|=|\alpha_1\dots\alpha_k|:=k.
$$

We denote monomials in $\Afr_n$ by
$$
\omega^\alpha:=\omega_{\alpha_1}\dots\omega_{\alpha_k},
\qquad \text{where $\alpha=\alpha_1\dots\alpha_k$.}
$$
In this notation,
 elements of $\Afr_n$ are finite sums:
\begin{equation}
	x=\sum_{\alpha\in\cA} c_\alpha \omega^\alpha.
	\label{eq:sum}
\end{equation}

We define  $\Z_+$-grading of $\Afr_n$ assuming $\deg \omega_j=1$.
So, $\deg \omega^\alpha=|\alpha|$.
The algebra $\Afr_n$ is a direct sum of homogeneous subspaces,
$$
\Afr_n=\bigoplus_{j=0}^\infty \Afr_n^{[j]}.
$$
Each element  $x\in\Afr_n$ is
a sum of its homogeneous components,
\begin{equation}
x=\sum_{j\ge 0} x^{[j]}, \qquad \text{where $x^{[j]}\in \Afr_n^{[j]}$}.
\end{equation}
We also denote by
$$
\phantom{A}^+\!\Afr_n:=\bigoplus_{j>0}^\infty \Afr_n^{[j]}
$$
the ideal of $\Afr_n$ consisting of elements with zero constant term.

By $\ov\Afr_n$ we denote the algebra of formal series 
in $\omega_j$, i.e., formal series of the form
\eqref{eq:sum}. We say that a sequence
$x_l=\sum_{\alpha\in\cA} c_\alpha^{(l)} \omega^\alpha\in \ov\Afr_n$
 {\it converges} to
$x=\sum_{\alpha\in\cA}c_\alpha \omega^\alpha\in \ov\Afr_n$ if for each $\alpha$ we have a convergence
of coefficients
$c_\alpha^{(l)}\to c_\alpha$.
So we get a topology of a direct product of a countable number
of copies
of  $\R$ enumerated by $\alpha\in\cA$.
 Thus, $\ov\Afr_n$ is a Fr\'echet  space%
 \footnote{Recall that a {\it Fr\'echet  space}
 is a complete topological vector space whose topology is defined by a countable family of seminorms, 
 in our case $\| x\|_{(\alpha)}=|c_\alpha|$, where $\alpha$ ranges in $\cA$.}.

We have
a well-defined {\it exponential} and a well-defined {\it logarithm}.
Namely, 
for $x$, $y\in \vphantom{A}^+\!\ov\Afr_n$
\begin{align*}
\exp x:&=\sum_{m=0}^\infty \frac{1}{m!}\, x^m, 
\qquad\text{notice, that we have $\exp(x)-1\in \vphantom{A}^+\!\ov\Afr_n$;}
\\
 \ln(1+y)&:=\sum_{m=1}^\infty \frac{(-1)^{m-1}}{m}\,y^m\,\,\in \vphantom{A}^+\!\ov\Afr_n,
\end{align*}
these series converge in $\ov\Afr_n$,
since in each homogeneous component $\Afr_n^{[j]}$
we actually have a finite sum. As usual,
$\ln(\exp(x))=x$, $\exp(\ln(1+y))=1+y$.

\sm

For $\xi\ge 0$ we define the endomorphism $\cM_\xi:\Afr_n\to \Afr_n$ by 
$$
\cM_\xi \omega_j=\xi \omega_j,
$$
or, equivalently,
$$
\cM_\xi \omega^\alpha=\xi^{|\alpha|} \omega^\alpha,
\quad \cM_\xi\bigl(\sum_j x^{[j]}\bigr)=\sum_j \xi^j x^{[j]}.
$$
Clearly, if $\xi>0$, then $\cM_\xi$ is an automorphism of $\Afr_n$.
It extends by continuity to an automorphism of 
 $\ov\Afr_n$. 

For $\xi=0$ we get the homomorphism
 $\ov\Afr_n\to \ov\Afr_n$ defined by  $\sum_j x^{[j]}\mapsto x^{[0]}$.

\sm

We also have a {\it canonical anti-automorphism} 
$$x\mapsto x^\circ$$
of $\Afr_n$ (and of $\ov\Afr_n$)
defined by
$$\omega^\circ= -\omega, \qquad
(xy)^\circ=y^\circ x^\circ,\qquad (x+y)=x^\circ+y^\circ
.$$
In other words,
$$
(\omega_{\alpha_1}\dots \omega_{\alpha_k})^\circ
=(-1)^k \omega_{\alpha_k}\dots \omega_{\alpha_1}.
$$

{\bf \punct  The coproduct.} Denote by $\Afr_n\otimes \Afr_n$
the   tensor product
of the algebra $\Afr_n$ with itself, i.e., the real associative algebra
with the basis $\omega^\alpha\otimes \omega^{\beta}$ and the multiplication
defined by
\begin{equation}
\omega^\alpha\otimes \omega^{\beta}\,\cdot\,
\omega^{\alpha'}\otimes \omega^{\beta'}=
\omega^\alpha \omega^{\alpha'}\otimes \omega^{\beta} \omega^{\beta'}.
\label{eq:tensor}
\end{equation}
The {\it coproduct} (see \cite{Reu}, Sect.~1.3-1.4)
$$\Delta:\,\Afr_n\to \Afr_n\otimes \Afr_n
$$
is a homomorphism of associative algebras
defined on generators by
$$
\Delta \omega_j=\omega_j\otimes 1+1\otimes \omega_j.
$$

For each $k\ge 0$ we denote by $\cP(k)$ the set of all subsets of $\{1,2,\dots,k\}$.
For any monomial $\omega^\alpha$ and any $I\in \cP(|\alpha|)$
we define  the monomial 
$$\omega^{\alpha|I}:=\omega_{\alpha_{i_1}}\dots \omega_{\alpha_{i_\mu}}$$
(i.e., we remove from $\omega_{\alpha_1}\dots\omega_{\alpha_k}$
symbols $\omega_{\alpha_j}$ with $j\notin I$).
Denote by $\backslash I$ the complement $\{1,2,\dots,k\}\setminus I$.





In this notation,
$$
\Delta(\omega^\alpha)=\sum_{I\subset \cP(k)}\omega^{\alpha|I}\otimes \omega^{\alpha|\backslash I}.
$$

The map $\Delta$ extends in the obvious way to
a continuous homomorphism 
$$\Delta:\,\ov\Afr_n\to \ov\Afr_n\otimes \ov\Afr_n.
$$
The symbol $\ov\Afr_n\otimes \ov\Afr_n$  here denotes 
 the algebra of formal series 
$$
\sum_{\alpha,\beta\in \cA}\kappa_{\alpha,\beta}
\,\omega^\alpha\otimes\omega^\beta
$$
with multiplication \eqref{eq:tensor}.

\sm

{\bf \punct Free Lie algebras.}
Denote by $\fr_n$ the free Lie algebra generated by $\omega_1$,
\dots, $\omega_n$ (on free Lie algebras, see a transparent introduction
in Serre \cite{Ser}, Ch.~IV, and a detailed exposition in  Reutenauer \cite{Reu}, see, also Bourbaki \cite{Bou}, Ch. 2). It is a positively graded Lie algebra
($\deg \omega_j=1$),
therefore we have a well-defined Lie commutator on the space
 $\ov\fr_n$ of formal series
$$
z=\sum_{j=1}^\infty z^{[j]}, \qquad \text{where $z^{[j]}\in \fr_n$, $\deg z^{[j]}=j$.}
$$

For any element of the Lie algebra $\fr_n$,
we assign  $z\in \Afr_n$ assuming that each
commutator $[a,b]$  in the expression for $z$
is the associative expression $ab-ba$ (see \cite{Reu}, Sect.~1.2). For instance,
\begin{align*}
&w_{17}\mapsto w_{17},\quad [\omega_1,\omega_2]\mapsto \omega_1\omega_2-\omega_2\omega_1,\\
&
\bigl[[\omega_1,\omega_2],\omega_4\bigr]\mapsto
(\omega_1\omega_2-\omega_2\omega_1)\omega_4
- \omega_4(\omega_1\omega_2-\omega_2\omega_1),
\\
& \bigl[[\omega_1,\omega_2], [\omega_7,\omega_8]\bigr]
\mapsto (\omega_1\omega_2-\omega_2\omega_1)
(\omega_7\omega_8-\omega_8\omega_7)
-(\omega_7\omega_8-\omega_8\omega_7)
(\omega_1\omega_2-\omega_2\omega_1),
\end{align*} 
etc.
We regard $\fr_n$  as
 the image of this map. So, 
 $$\fr_n\subset \Afr_n,\quad 
 \ov\fr_n\subset \ov\Afr_n.$$

 The following statement is well-known, see, for example, 
 \cite{Ser}, Sect.~IV.4, \cite{Reu}, Theorem 1.4 (the statement c)
 was discovered by Friedrichs \cite{Fri}, p.~203, see also \cite{Mag}).
 
 \begin{theorem}
 \label{th:friedrichs}
{\rm a)} The map $\fr_n\to \Afr_n$ is an embedding.

\sm

{\rm b)} The algebra $\Afr_n$ is the universal enveloping algebra of
$\fr_n$.

\sm

{\rm c)} An element $z\in \Afr_n$ is contained in $\fr_n$
if and only if
$$
\Delta(z)=z\otimes 1+1\otimes z.
$$ 

{\rm d)} The  statements {\rm a)}, {\rm c)} take place for
$\ov\fr_n$ and $\ov \Afr_n$.
 \end{theorem}
 
{\sc Remark.} Let $\frg$ be a Lie algebra, $U(\frg)$ be its
universal enveloping algebra. The diagonal embedding $\frg\to\frg\oplus\frg$ induces a homomorphism 
$U(\frg)\to U(\frg)\otimes U(\frg)$,
we have $x\mapsto x\otimes 1+1\otimes x$ for any  $x\in \frg$,
see Dixmier \cite{Dix}, Sect.~2.7 and 2.8.16-17. 
If $\frg=\fr_n$ this gives the coproduct $\Delta$. 
\hfill $\boxtimes$

\sm
 
{\bf \punct Free Lie groups.} 
Recall Campbell--Hausdorff theorem, see, e.g., 
\cite{Ser}, Sect.~IV.7, \cite{Reu}, Ch.~3.

\begin{theorem}
For any $x$, $y\in \ov\fr_n$, there exists a unique
$z\in \ov\fr_n$ such that
$$
\exp(x)\exp(y)=\exp(z).
$$
\end{theorem}

Therefore, the image of $\ov\fr_n$ under
the exponential map $\ov\fr_n\to \ov\Afr_n$ is a group with respect to the multiplication. Denote it by $\ov\Fr_n$.

Clearly, $x\in\ov\fr_n$ satisfies $x^\circ=-x$.
Therefore,
$$\exp(x)^\circ=\exp(x^\circ)=\exp(-x)=\exp(x)^{-1}.$$
Thus, for $g\in \ov\Fr_n$ we have 
\begin{equation}
g^{-1}=g^\circ.
\label{eq:g-circ}
\end{equation}

 Elements of this group can be characterized 
in the following more convenient  way (see Ree \cite{Ree}, see also \cite{Reu}, Theorem 3.2):

\begin{theorem}
An element $g\in\ov\Afr_n$ is contained in
$\ov \Fr_n$ if and only if
\begin{equation}
\Delta g=g\otimes g.
\label{eq:ggg}
\end{equation}
\end{theorem}

Let $g=\sum c_\alpha \omega^\alpha$. Then the condition \eqref{eq:ggg} can be 
written as a system of quadratic equations
\begin{equation}
c_\alpha c_\beta=\sum_I c_{\mu(I)},
\label{eq:system}
\end{equation}
where the summation is taken over all 
 $|\alpha|$-element subsets $I\subset\cP(|\alpha|+|\beta|)$
and $\mu(I)\in \cA$  is defined from the conditions:

\sm

--- $|\mu|=|\alpha|+|\beta|$;

\sm

--- $\omega^{\mu(I)|I}=\omega^\alpha$, $\omega^{\mu(I)|\backslash I}=\omega^\beta$.

\sm

So, the summation is taking over all the {\it shuffles} of monomials
$\omega^\alpha$, $\omega^\beta$.

\begin{corollary}
The group $\ov\Fr_n$ is closed in $\ov\Afr_n$.
\end{corollary}

{\sc Remark.} 
Apparently, firstly the free Lie group $\Fr_n$ was introduced by Malcev \cite{Mal} as a special case of his Lie-zation of generalized nilpotent
groups.
On some related groups, see \cite{Ner}.
\hfill $\boxtimes$

\sm

{\bf \punct Differentiations of paths and ordered exponentials%
\label{ss:ordered}.} Below we need the following simple arguments
(see, e.g., \cite{Mag}), which allow to apply  Volterra type integral equations to the free Lie groups.  We must fix a  degree of generality,
which is necessary below for the proof of Proposition 
\ref{pro:image-dense},
 so we present some obvious details.

\begin{proposition}
\label{pr:derivative}
 Let $q(t)=\sum_j q^{[j]}$ be a path in $\ov\Fr_n$, let $q(t)$ has
 a derivative at $t=t_0$, i.e., each homogeneous  component $q^{[j]}(t)$
 has  a derivative.
Then 
\begin{equation}
q(t_0)^{-1}q'(t_0)\in \ov\fr_n.
\label{eq:qq-prime}
\end{equation}
\end{proposition}

{\sc Proof.}
Without loss of generality, we can assume $t_0=0$.
Denote $q:=q(0)$, $q':=q'(0)$.
Then 
\begin{multline*}
\Delta q(\epsilon)=\Delta \bigl(q +\epsilon q'+o(\epsilon)\bigr)=
\Delta q \cdot \bigl(1+\epsilon \Delta (q^{-1} q')+o(\epsilon)\bigr)
=\\= (q\otimes q)\cdot \bigl(1+\epsilon \Delta (q^{-1} q')+o(\epsilon)\bigr), 
\end{multline*}
where $o(\epsilon)$ denotes an expression, which is
 $o(\epsilon)$ in each homogeneous
component.  Actually, in $j$-th homogeneous component
we have in the rind hand side  a polynomial
in $g^{[1]}$, \dots $g^{[j]}$, $(g')^{[1]}$, \dots, $(g')^{[j]}$
plus $o(\epsilon)$.
On the other hand,
\begin{multline*}
q(\epsilon)\otimes q(\epsilon)=(q\otimes q)
 \cdot\bigl(1+\epsilon q^{-1} q'+o(\epsilon) )\otimes (1+\epsilon q^{-1} q'+o(\epsilon) \bigr)
 = \\=
 (q\otimes q)\cdot \bigl(1+\epsilon (1\otimes q^{-1} q'
 + q^{-1} q'\otimes 1)+o(\epsilon)\bigr).
\end{multline*}
Comparing coefficients at $\epsilon$
and keeping in mind Friedrichs theorem \ref{th:friedrichs}, we come to the desired statement.
\hfill $\square$

\begin{proposition}
\label{pr:ordered}
Let $\gamma(t)=\sum \gamma^{[j]}(t)$ be a measurable map
 $[0,1]\to\ov\fr_n$. Let all components
$\gamma^{[j]}(t)$ be Lebesgue integrable $\Afr^{[j]}_n$-valued functions.

\sm

{\rm a)}
Then there exists 
a unique   solution  (the '{\rm ordered exponential}')
$$
\ex[\gamma;t]=1+\sum_{j\ge 1}\ex^{[j]}[\gamma;t] \in \ov\Fr_n
$$  
of the Cauchy problem
 \begin{align}
\frac d{dt}
\ex[\gamma;t]&= \ex[\gamma;t]\, \gamma(t);
\label{eq:cauchy}
\\ \ex[\gamma;0]&=1.
\label{eq:cauchy2}
\end{align}
such that all $\ex^{[j]}[\gamma;t]$ are absolutely continuous
in $t$,
and \eqref{eq:cauchy} holds almost everywhere.

\sm

{\rm b)} The ordered differential $\ex[\gamma;t]$
is a unique continuous solution of the
 following Volterra integral equation:
\begin{equation}
\ex[\gamma;t]=1+\int_0^t \ex[\gamma;\tau]\,\gamma(\tau)\,d\tau.
\label{eq:volterra}
\end{equation}

\end{proposition}

{\sc Proof.} a)
We must solve the equation
$$
\frac d{dt}
\ex[\gamma;t]=
\bigl(1+\ex[\gamma;t]^{[1]}+\ex[\gamma;t]^{[2]}+\dots\bigr)\,
\bigl( \gamma(t)^{[1]}+\gamma(t)^{[2]}+\dots\bigr).
$$
So, for $j>0$ we have
\begin{align*}
\frac d{dt}\ex[\gamma;t]^{[j]}&=
\gamma(t)^{[j]}+ \gamma(t)^{[j-1]}  \ex[\gamma;t]^{[1]}+
\dots + \gamma(t)^{[1]}  \ex[\gamma;t]^{[j-1]};
\\ 
 \frac d{dt}\ex[\gamma;0]^{[j]}&=0,
\end{align*}
notice that $\ex[\gamma;t]^{[j]}$ is absent in the right hand side
of the differential equation.
 We  find $\ex[\gamma;t]^{[1]}$,
$\ex[\gamma;t]^{[2]}$, $\dots$ recurrently by indefinite integrations;
recall that a primitive of a Lebesgue integrable function
is unique (up to an additive constant) in the class of absolutely continuous functions.

It is easy to see that if $\gamma_k(t)\in \ov\fr_n$ is a sequence
such that all its homogeneous components $\gamma_k^{[j]}(t)$ converge 
to $\gamma^{[j]}(t)$ in the $L^1$-sense, then $\ex[\gamma_k;t]^{[j]}$ converge uniformly  to $\ex[\gamma;t]^{[j]}$.
 
To verify
$\ex[\gamma;t]\in \ov\Fr_n$,
we take a sequence of piecewise constant functions $\gamma_k(t)$
 convergent  to  $\gamma(t)$ in $L^1$-sense
 and pass to the limit as $k\to\infty$. 
 

\sm

b) Taking integral $\int_0^t$ of  both sides of differential equation
\eqref{eq:cauchy} and keeping in mind absolute continuity
of $\cE[\dots]$, we come
to \eqref{eq:volterra}. On the other hand,
if in the integral equation \eqref{eq:volterra}
a function $\gamma(t)$ is in $L^1$ and $\ex$ is continuous,
then $\int_0^t(\dots)$ is absolutely continuous; differentiating
both sides by $t$, we get the differential equation \eqref{eq:cauchy}.
\hfill $\square$

\section{The groups $\Fr_n^\infty$ and $\Fr_n^\xi$}

\COUNTERS

{\bf \punct Algebras $\Afr^\xi_n$ and $\Afr^\infty_n$.}
Fix $\xi>0$. We define a norm $\|x\|_\xi$ of
 $$x=\sum_\alpha c_\alpha \omega^\alpha\in \ov\Afr_n$$
by
$$\|x\|_\xi=\sum_\alpha \xi^{|\alpha|}|c_\alpha|.$$
Denote by $\Afr_n^\xi$ the Banach space of all $x$ such that
$\|x\|_\xi<\infty$. We 
have
$$
\|x+y\|_\xi\le \|x\|_\xi+\|y\|_\xi,
\qquad \|xy\|_\xi\le \|x\|_\xi\,\|y\|_\xi, 
\qquad \|\,[x,y]\,\|_\xi\le 2\,\|x\|_\xi\,\|y\|_\xi.
$$
The second inequality is clear for $\xi=1$,
for general $\xi$ it follows from
$$
\|x\|_\xi=\|\cM_{\xi} x\|_1.
$$

Therefore, {\it  $\Afr_n^\xi$ are Banach algebras}.
If $\eta>\xi$, then $\Afr_n^\eta\subset \Afr_n^\xi$,  

We define the  algebra $\Afr_n^\infty$
by
$$
\Afr_n^\infty:=\bigcap_{\xi>0} \Afr_n^\xi.
$$
We equip it with the topology of a polynormed space and get a Fr\'echet space. Clearly, the multiplication in $\Afr_n^\infty$  is jointly continuous.

\sm

For $0<t<1$, the map $\cM_t$ is an endomorphism of $\Afr_n^\xi$,
for all $t>0$ it is an automorphism of $\Afr_n^\infty$.
The map $x\mapsto x^\circ$ is an anti-automorphism.

\sm

For  a finite-dimensional normed space $V$,
denote by $\Mat(V)$
 the Banach algebra of linear operators in $V$.

\begin{proposition} 
\label{pr:map}
Let $B_1$, \dots, $B_n$ be  operators in $V$.
 
 \sm
 
{\rm a)}   The homomorphism $\Afr_n\to \Mat(V)$ 
 defined by $\omega_j\mapsto B_j$ extends by continuity
 to a homomorphism $\Afr_n^\infty\to\Mat(V)$.
 
\sm 
 
 {\rm b)} If  norms $\tri B_j\tri$ of all $B_j$ are  $\le\xi$,
 then this homomorphism  extends by continuity
 to a homomorphism $\Afr_n^\xi\to \Mat(V)$.
\end{proposition}

{\sc Proof.} a) follows from b). The statement b) is obvious:
for $x=\sum_\alpha c_\alpha \omega^\alpha\in \Afr_n^\xi$,
the series 
$$
\sum_{\alpha\in \cA} c_\alpha B^\alpha=\sum_{\alpha\in \cA}
c_\alpha B_{\alpha_1}\dots B_{\alpha_k} 
$$
absolutely converges, since the series of norms is
$$
\sum_{\alpha\in \cA}
|c_\alpha|\, \tri B_{\alpha_1}\dots  B_{\alpha_k}\tri
\le
\sum_{\alpha\in \cA}
|c_\alpha|\, \tri B_{\alpha_1}\tri\dots \tri B_{\alpha_k}\tri
\le \sum_{\alpha\in \cA}
|c_\alpha|\, \xi^{|\alpha|}=\|x\|_\xi.
$$

{\sc Remark.} Completions of universal enveloping algebras were considered
by different authors for different reasons. See Hohschild \cite{Hoch}
(see, also, \cite{Dix}, 2.8.16-17),
Rashevskii \cite{Rash}, 
Litvinov \cite{Lit}, Goodman \cite{Goo1}--\cite{Goo3},
Taylor \cite{Tay},  Zhelobenko \cite{Zhe}, Pirkovskii \cite{Pir},
Neretin \cite{Ner}.  A nilpotent Lie group $G$
can be embedded to an appropriate completion of the corresponding 
enveloping algebra $U(\frg)$, see Goodman \cite{Goo1}--\cite{Goo3}. 
 We follow \cite{Ner}, where the main topic were completions of braid groups. See, also, \cite{Ner2}.
\hfill $\boxtimes$

\begin{lemma}
\label{l:continuity}
The maps $[0,1]\times\Afr_n^\xi\to \Afr_n^\xi$,
$[0,1]\times\Afr_n^\infty\to \Afr_n^\infty$
given by
$
(t, x)\to \cM_t x
$
are continuous.
\end{lemma}

{\sc Proof.} It is sufficient to verify the statement for $\xi=1$.
Let $t_k\to t$, $x_k\to x$. Then
\begin{multline*}
\|\cM_{t_k}x_k- \cM_{t}x \|_1\le \|\cM_{t_k}x_k- \cM_{t_k}x \|_1
+
\|\cM_{t_k}x- \cM_{t}x \|_1
=\\=\sum_{j\ge 0} \|t_k^j x_k^{[j]}- t_k^j x^{[j]}\|_1
+ \sum_{j\ge 0} \|t_k^j x^{[j]}-t^j x^{[j]}\|_1
=\\= \sum_{j\ge 0} t_k^j\, \|x_k^{[j]}- x^{[j]}\|_1
+ \sum_{j\ge 0} |t_k^j-t^j|\cdot  \|x^{[j]}\|_1.
\end{multline*}
The first series is $\le\sum  \|x_k^{[j]}- x^{[j]}\|_1=\|x_k-x\|_1\to 0$.
The second series is dominated by a convergent series 
$\sum \|x^{[j]}\|_1=\|x\|_1$. By the dominated convergence
theorem, the termwise convergence of summands
to 0 implies a convergence of their sums to 0.
\hfill $\square$

\sm 

{\bf\punct Completed free Lie algebras and the corresponding Lie groups.}
We define
$$
\fr_n^\xi:=\ov\fr_n \cap \Afr_n^\xi,\quad 
\fr_n^\infty:=\ov\fr_n \cap \Afr_n^\infty,
$$  
clearly, these subspaces are Lie algebras.  
The algebras  $\fr_n^\xi$ are Banach spaces, the algebra $\fr_n^\infty$
is a Fr\'echet space. In both cases, 
the brackets $[\cdot,\cdot]$ are jointly continuous.

Following \cite{Ner}, we also define  groups
$$
\Fr^\xi_n:=\ov \Fr_n\cap \Afr_n^\xi, \qquad \Fr^\infty_n:=\ov \Fr_n\cap \Afr_n^\infty=\cap_{\xi>0} \Fr_n^\xi,
$$
by \eqref{eq:g-circ} these intersections are closed with respect
the inversion $g\mapsto g^{-1}$. We equip these groups with
 topologies induced from $\Afr_n^\xi$ and $\Afr_n^\infty$ respectively.

Clearly, all groups $\Fr^\xi_n$, where $\xi>0$, are canonically isomorphic. 

\sm

\begin{proposition}
The subsets  $\Fr^\xi_n\subset \Afr^\xi_n$,
 $\Fr^\infty_n\subset \Afr^\infty_n$ are closed.
 Groups $\Fr^\xi_n$, $\Fr^\infty_n$ are Polish%
 \footnote{A Polish space is a topological space 
 homeomorphic to a complete separable metric space. A Polish
 group is a topological group, which is a Polish space.}.
\end{proposition}

{\sc Proof.}
The condition $\Delta g=g\otimes g$ is a countable system of 
conditions \eqref{eq:system}. Each condition determines a quadratic surface
in a finite dimensional space 
$$
\bigl(\Afr_n^{[|\alpha|]}+\Afr_n^{[|\beta|]}\bigr)\oplus
\Afr_n^{[|\alpha|+|\beta|]},
  $$ 
  which is a quotient of the whole space, and therefore
  each condition \eqref{eq:system} determines a closed subset.
  
  A closed subset in a Polish space is Polish.
\hfill $\square$

 \begin{proposition}
 The groups $\Fr^\xi_n$ and $\Fr^\infty_n$ are
 contractible {\rm(}and therefore linear connected and simply connected{\rm)}.
 \end{proposition}
  
{\sc Proof.} If $g$ satisfies the condition $\Delta g=g\otimes g$,
then $\cM_t g$ satisfies the same condition. 
By Lemma \ref{l:continuity},  the maps $[0,1]\times \Fr_n^\xi\to \Fr_n^\xi$, $[0,1]\times \Fr_n^\infty\to \Fr_n^\infty$
given by $(t,g)\mapsto \cM_t g$ are continuous.
For $t=1$ we have an identical map from the group to itself, for $t=0$
the group is mapped to the unit.  
%
\hfill $\square$

\sm 

Since $\Afr_n^\xi$ are Banach algebras, the exponential
$x\mapsto \exp x$ is a well-defined map from
$\vphantom{A}^+\!\Afr_n^\xi$ to the affine subspace 
$1+\vphantom{A}^+\!\Afr_n^\xi$. Since this is valid 
for all $\xi$, we have a map 
$\vphantom{A}^+\!\Afr_n^\infty\to 1+\vphantom{A}^+\!\Afr_n^\infty$.
Restricting these maps to the free Lie algebras, we get well-defined exponentials
$$
\exp:\,\fr_n^\xi\to \Fr_n^\xi, \qquad 
\fr_n^\infty\to \Fr_n^\infty.
$$

Similarly, 
the logarithm $y\mapsto \ln(1+y)$
 is a  well-defined map from the unit ball $\|y\|_\xi<1$ in  
 $\vphantom{A}^+\!\Afr_n^\xi$ to $\vphantom{A}^+\!\Afr_n^\xi$.
 
  If also $1+y\in \Fr_n^\xi$,
   then we have $\ln(1+y)\in \ov\fr_n$,
 and so $\ln (1+y)\in \fr_n^\xi$.

\begin{corollary}
\label{cor:generated}
The group $\Fr_n^\xi$ is generated by the image
of the exponential map $\exp:\fr_n^\xi\to \Fr_n^\xi$.
\end{corollary}

Indeed, the image of exponential map contains a neighborhood of unit
in $\Fr_n^\xi$,
and the group is connected.
\hfill $\square$

\sm

{\sc Remark.} So, {\it $\Fr_n^\xi$ is a Banach--Lie group}. I.e.,
it is 
a topological group, whose neighborhood of the unit is modeled by an open
domain in a Banach space, and the multiplication is analytic.
Many elementary statements about  Lie groups and Lie algebras
 survive for Banach--Lie groups, see Bourbaki \cite{Bou}, Ch. 3,
 see also \cite{Har}, \cite{Swi}. 
However, numerous
'infinite-dimensional Lie groups', which were subjects and tools  in mathematics and mathematical physics
during the last 60 years, usually are not contained in this class.
Note, that our Banach--Lie group $\Fr_n^\xi$ is an element of a big family
of related groups, see \cite{Ner}, Sect. 4, \cite{Ner2}. 
I do not know if $\Fr_n^\infty$ is an infinite-dimensional manifold
in some formal sense.
\hfill $\boxtimes$

\begin{observation} The image of the exponential
map $\fr_n^\infty\to \Fr_n^\infty$ does not contain 
a neighborhood of the unit.
\end{observation}

{\sc Proof.}  The  curve 
$$\mu(t)=\exp(t \omega_1)\,\exp(t\omega_2)$$
 is continuous
 in $\Fr_n^\infty$ and passes through the unit ($\mu(0)=1$).
Let us show that $\mu(t)$ is not contained in the image of the exponential
map for all $t>0$.  
  Consider its logarithm,
 \begin{multline*}
 \ln\mu(t)=
 \ln \Bigl[ \Bigl(1+t\omega_1+\sum_{i\ge 2}
  \frac{t^i}{i!}\,\omega_1^i\Bigr)
  \Bigl(1+t\omega_2+\sum_{j\ge 2} \frac{t^j}{j!}\,\omega_2^j\Bigr)\Bigr]
  =\\= \ln\bigl[1+t(\omega_1+\omega_2)+\dots\bigr]=
  \sum_{k=1}^\infty \frac{(-1)^{k-1}\, t^k}{k}\,(\omega_1+\omega_2)^k+
  \dots
 \end{multline*}
 Taking $k=2m$, we observe summands 
 $\Bigl(-\frac{t^{2m}}{2m}(\omega_1\omega_2)^m\Bigr)$ 
 and $\Bigl(-\frac{t^{2m}}{2m}(\omega_2\omega_1)^m\Bigr)$. 
 The whole series for $\ln\mu(t)$ does not contain other similar 
 summands with $(\omega_1\omega_2)^m$, $(\omega_2\omega_1)^m$
 --- all other monomials of the expansion contain subwords
 of the form $\omega_1^i$ with $i\ge 2$
 or  $\omega_2^j$  with $j\ge 2$. So the series
 $ \ln\mu(t)$
 contains a subseries
 $$
 -\sum_{m>0}\frac{t^{2m}}{2m}\Bigl((\omega_1\omega_2)^m+(\omega_2\omega_1)^m\Bigr).
 $$
 If $t>\xi^{-1}$, then $\ln\mu(t)\notin \Afr_n^\xi$.
 Therefore, $\ln\mu(t)\notin \fr_n^\infty$ for all $t>0$.
 \hfill $\square$
 
  
\begin{proposition}
\label{pro:image-dense}
The subgroup in $\Fr_n^\infty$ generated by the image
of the exponential map $\fr_n^\infty\to  \Fr_n^\infty$ is 
dense in $\Fr_n^\infty$.
\end{proposition}
 
The statement is proved in the next subsection.

\sm

{\bf \punct Ordered exponentials in $\Fr_n^\xi$ and $\Fr_n^\infty$.}
Here we prove Proposition \ref{pro:image-dense}.

\begin{lemma}
\label{l:dependence}
Let $\gamma:[0,1]\to \fr_n^\xi$ be a measurable map.
Denote 
$$I_\xi(\gamma):=\int_0^1 \|\gamma(t)\|_\xi\,dt.$$

{\rm a)} 
Let $I_\xi(\gamma)<\infty$.
Then its ordered exponential  $\ex[\gamma;t]$ is  a continuous map
$[0,1]\to \Fr_n^\xi$.

\sm 

{\rm b)} Let $\gamma$, $\gamma_j: [0,1]\to \fr_n^\xi$ be  measurable maps,
such that $I_\xi(\gamma)$, $I_\xi(\gamma_j)<\infty$.
Let $I_\xi(\gamma_j-\gamma)\to 0$.
Then  $\|\ex[\gamma_k;t]- \ex[\gamma;t]\|_\xi$ converges to $0$ uniformly in $t$ as $k\to\infty$.  
\end{lemma}

{\sc Proof.}
 a) We  consider 
the Banach space $W_\xi=C([0,1], \Afr_n^\xi)$ of continuous functions
$[0,1]\to \Afr_n^\xi$
equipped 
with the uniform norm $\max_t \|f(t)\|_\xi$.
Consider a Volterra integral operator $A_\gamma$ in this space given by
$$
A_\gamma f(t)=\int_0^t f(\tau)\,\gamma(\tau)\, d\tau,
$$
its norm is $\le I_\xi(\gamma)$. Indeed,
\begin{multline*}
\|A_\gamma f\|_{W_\xi}=\max\limits_{t\in[0,1]}
\|\int_0^t f(\tau)\,\gamma(\tau)\, d\tau\|_\xi
\le  \max\limits_{t\in[0,1]}\int_0^t \|f(\tau)\|_\xi\,\|\gamma(\tau)\|_\xi\, d\tau=
\\=\int_0^1 \|f(\tau)\|_\xi\,\|\gamma(\tau)\|_\xi\, d\tau
\le \int_0^1 \|\gamma(\tau)\|_\xi\, d\tau  \cdot \max\limits_{\tau\in[0,1]}\|f(\tau)\|_\xi
=I_\xi(\gamma)\cdot\|f\|_{W_\xi}.
\end{multline*}
The same estimate gives
$$
\|A_\gamma-A_\mu\|_{W_\xi}\le I_\xi(\gamma-\mu).
$$
To find the ordered exponential $f=\ex[\gamma;t]$, we must solve the integral equation $f=1+A_\gamma f$, see \eqref{eq:volterra}.
Spectral radii of Volterra operators in Banach spaces
of continuous functions are $0$, see textbooks on functional
analysis, e.g. \cite{KA}, Subsect. XIII.6.2,
the same estimates hold for spaces of continuous functions with 
values in Banach spaces.
So $(1-A_\gamma)^{-1}$ is a bounded
operator. Therefore,
$$
\ex[\gamma;t]= (1-A_\gamma)^{-1}\cdot 1
\in C([0,1], \Afr_n^\xi).
$$

\sm 

b) The inversion of a linear operator $S\mapsto S^{-1}$ in a Banach space is 
a continuous operation with respect to the norm topology, see, e.g., \cite{KA}, Subsect. V.4.6. So, the sequence
$(1-A_{\gamma_k})^{-1}\cdot 1$ uniformly converges to $(1-A_{\gamma})^{-1}\cdot 1$.
 \hfill $\square$


 

\begin{lemma}
\label{l:ordered}
For $g\in \Fr_n^\infty$, we have 
$$\gamma(t):=(\cM_t g)^{-1} \cdot\frac{d}{dt} \cM_t g\,\,
\in \fr_n^\infty,
$$
$\gamma(t)$ is continuous,
and
$$
\ex[\gamma;t]=\cM_t g.
$$
\end{lemma}  
  
{\sc Proof.} Clearly, the derivative 
$$ \frac{d}{dt} \cM_t g= \frac{d}{dt} \sum_j t^j g^{[j]}
=\sum_j j t^{j-1} g^{[j]}$$
is contained in $\Afr^\infty_n$ and is continuous in $t$.
By Proposition \ref{pr:derivative}, we get the first claim.
On the other hand, maps $t\mapsto \cM_tg^{[j]}=t^j g^{[j]}$
are absolutely continuous on finite intervals, and we refer to Proposition \ref{pr:ordered}.
 \hfill $\square$

\sm


{\sc Proof of Proposition \ref{pro:image-dense}.} 
Denote  by $\frF\subset \Fr_n^\infty$ the subgroup
 generated by the image
of the exponential map.

 If $\gamma$ is constant, when $\ex[\gamma;t]=\exp(\gamma t)\in \frF$. If $\gamma$ is locally constant, then $\ex[\gamma;t]$
 is a product of several exponentials and also is contained in $\frF$.
 
 Fix $g\in\Fr_n^\infty$, denote  $\gamma(t):=(\cM_t g)^{-1} (\cM_t g)'\in 
 \fr_n^\infty$ as in Lemma \ref{l:ordered}. So, $\ex[\gamma;1]=g$.
  Fix an increasing sequence
 $\xi_k\to\infty$ and a decreasing sequence $\epsilon_k\to 0$.
 For each $k$ we choose a piece-wise constant function $\gamma_k$
 such that 
 $$\sup_{t\in[0,1]}\|\gamma_k(t)-\gamma(t)\|_{\xi_k}<\epsilon_k.$$ 
 To do this, for each $N$ we 
 set $\Gamma_{N}(t):=\gamma(p/N)$ if $p/N\le t<(p+1)/N$.
 Since the function $\gamma(t)$ is continuous 
as a function $[0,1]\mapsto \Afr_n^{\xi_k}$, it is uniformly
continuous, and therefore 
$\sup_{t\in[0,1]}\|\gamma(t)-\Gamma_N(t)\|_{\xi_k}$
tends to 0 as $N\to\infty$. 
  For sufficiently large $L=L(k)$
we have $\sup_{t\in[0,1]}\|\gamma(t)-\Gamma_L(t)\|_{\xi_k}<\epsilon_k$,  
   and we can take 
 $\gamma_{k}(t):=\Gamma_{L(k)}(t)$.

 Since the norms $\|\cdot\|_{\xi_k}$ are increasing, 
  for a fixed $\xi_m$ we have 
  $$\| \gamma_k(t)-\gamma(t)\|_{\xi_m}
\le \| \gamma_k(t)-\gamma(t)  \|_{\xi_k}  
  \to 0$$
   as
  $k\to \infty$. By Lemma \ref{l:dependence}, the sequence
$\ex[\gamma_k;1]\in \frF$ converges to $\ex[\gamma;1]=g$
in each $\Fr_n^\xi$ and so converges in $\Fr_n^\infty$.
\hfill $\square$

\sm  
  
{\bf \punct Universality.} Let $G$ be a connected 
 real Lie group, %
let $\frg$ be its Lie algebra, let $\Exp:\frg\to G$ be the exponential
map. Recall that its image contains a neighborhood
of the unit in $G$. Let $X_1$, \dots, $X_n\in \frg$.
Then we have a canonical homomorphism $\pi:\fr_n\to\frg$
such that $\pi(\omega_j)=X_j$.   
  
 \begin{theorem}
 \label{th0}
 There exists a unique homomorphism 
$\Fr_n^\infty\to G$ satisfying 
\begin{equation}
\Pi(\exp(z))=\Exp(\pi(z))
\qquad\text{for any $z\in \fr_n$.}
\label{eq:exp-exp}
\end{equation}
If $X_1$, \dots, $X_n$ generate $\frg$, then
$\Pi$ is surjective.
 \end{theorem}
  
  We also formulate similar statements for the groups
  $\Fr_n^\xi$.
  
\begin{theorem}
\label{th}
Under the same conditions:

\sm 

{\rm a)} 
There is $\xi>0$ such that there exists a 
homomorphism
$\Pi:\Fr_n^\xi\to G$ satisfying the condition
\eqref{eq:exp-exp}.
If $X_1$, \dots, $X_n$ generate $\frg$, then
$\Pi$ is surjective. 

\sm 

{\rm b)} If for some $\xi>0$ there exists a homomorphism
$\Pi:\Fr_n^\xi\to G$ satisfying \eqref{eq:exp-exp},
then such homomorphism is unique.
\end{theorem}  

{\sc Remarks.}
a)
These statements are very simple but I have not seen
them in  literature. Pestov \cite{Pes} constructed 
a separable Banach--Lie
group, such that any separable Banach--Lie group can be realized
as its quotient group. Formally, our statements do not follow
from Theorem 3.5 of \cite{Pes}.  In any case, we provide an explicit construction
of the universal group, 
 and
our arguments  are very elementary. 

\sm

b) In particular, we see that unitary groups $\mathrm{U}(k)$ and
orthogonal groups $\SO(k)$ (having nontrivial topology) are quotient groups of our contractible 
groups $\Fr_n^\xi$ and $\Fr_n^\infty$ (for any $n\ge 2$). 
This is unusual from the point of view of classical theory
of Lie groups.
May be here a correct classical analogy are classifying spaces,
which are 
quotients  of contractible spaces by group actions.
\hfill $\boxtimes$  

\begin{lemma}
\label{l:th}
The statements of Theorems {\rm \ref{th0}--\ref{th}} hold under the following additional condition:
the group $G$ is linear, i.e., admits a faithful finite-dimensional representation.
\end{lemma}

{\sc Proof of Lemma \ref{l:th}.} Let $f$ be a map $A\to B$, let $C\subset A$.
In this proof, we denote by $f(C)$ the image of $C$ under the map $f$.

\sm

 Consider a faithful linear representation
$R$
of the group $G$ in a fi\-ni\-te-di\-men\-si\-o\-nal real linear
space $V$. Denote the corresponding representation of $\frg$
by $\rho$. Denote by $\GL(V)$ the group of invertible
linear operators in $V$, by $\frg\frl(V)$  the Lie algebra of linear operators in $V$. 
 Consider the group $R(G)\subset \GL(V)$
and the Lie algebra
$\rho(\frg)\subset \frg\frl(V)$.
  For any $Y\in \frg$,
we have 
$$
R(\Exp(Y))=\exp(\rho(Y)).
$$

It may happened that the subgroup $R(G)$ 
is not closed in $\GL(V)$. Formally we must distinguish
the group $R(G)$ with the topology arising from the isomorphism
$G\to R(G)$ and the group $R(G)^{twin}$ with the topology induced
from $\GL(V)$. In this case, a small neighborhood of the unit in
$R(G)^{twin}$ is disconnected, the connected component of the unit
is a neighborhood of unit in $R(G)$. Therefore, for any
 linear connected group $\Gamma$ a continuous homomorphism $\Gamma\to R(G)^{twin}$
is continuous as a homomorphism $\Gamma\to R(G)$.
Now we are ready to start the proof.

\sm 

Let us verify the statement a) of Theorem \ref{th}.  Consider some norm in $V$
and denote by $\tri\cdot\tri$ the corresponding
operator norm.
 Choose $\xi\ge\max \tri\rho(X_j)\tri$.
By Proposition \ref{pr:map}, the map $\omega_j \mapsto \rho (X_j)$
determines a homomorphism of Banach algebras
 $\Theta:\Afr_n^\xi\to \Mat(V)$.
Restricting $\Theta$ to the Lie algebra 
$\fr_n^\xi$, we get a continuous  homomorphism from
$\fr_n^\xi$ to the Lie algebra $\frg\frl(V)$. 
It coincides with $\rho\circ\pi$ on $\fr_n$.

Restricting $\Theta$ to $\Fr_n^\xi$, we get a homomorphism
 $\Fr_n^\xi\to \GL(V)$. For elements
$\exp(z)$, where $z\in \fr_n$,
we have 
$$\Theta(\exp(z))=\exp(\Theta(z))\in R(G).$$
By the continuity, this is so for $z\in\fr^\xi_n$.
 By Corollary \ref{cor:generated},
this implies $\Theta(\Fr_n^\xi)\subset R(G)$.
So, $\Theta$ determines a continuous homomorphism $\Fr_n^\xi\to R(G)^{twin}$.
Therefore, it is continuous as a homomorphism  $\Fr_n^\xi\to R(G)$.

We define $\Pi:\Fr^\xi_n\to G$ from the condition:
$$
\Theta=R\circ \Pi,
$$
keeping in mind that $R$ is an isomorphism between $G$ and $R(G)$.

\sm

Next, let $X_1$, \dots, $X_n$ generate $\frg$. Equivalently,
let $\pi:\fr_n\to\frg$ be surjective.
So  $\rho\circ\pi(\fr_n)=\rho(\frg)$, by the continuity we get
 $$\Theta(\fr_n^\xi)=\rho(\frg).$$ 
Therefore, the set $\Theta(\exp(\fr_n))\subset R(G)$ contains a neighborhood of
the unit. Since $G$ is connected, the map $\Theta:\Fr_n^\xi\to R(G)$
is surjective. 

\sm

Next, let us verify the statement b). The set $\exp(\fr_n)$ is dense in $\exp(\fr_n^\xi)\subset \Fr_n^\xi$.
Therefore, $\Pi$ is uniquely determined in a neighborhood
of the unit in $\Fr_n^\xi$. Since $\Fr_n^\xi$ is connected, $\Theta$
is uniquely determined on the whole $\Fr_n^\xi$.


\sm 

It remains to verify the statement of Theorem \ref{th0}. Considering the homomorphism $\Pi:\Fr_n^\xi\to G$ 
constructed above and restricting it to $\Fr^\infty_n$
we get existence.

 The uniqueness follows from Proposition \ref{pro:image-dense}.

If $\pi$ is surjective, then $\Pi(\exp(\fr_n))=\Pi(\frg)$
contains a neighbourhood of the unit. Hence, the image of $\Pi$
coincides with $G$.
\hfill $\square$

 
\sm

{\sc Proof of Theorems \ref{th0}--\ref{th}.}
%
By  Ado theorem, any finite-dimensional Lie algebra admits a faithful
representation. 
 Equivalently, for any connected Lie group $G$ some quotient
 of $G$ by a discrete subgroup $Z$ in its center is a linear group.
%
By Lemma \ref{l:th}, we have a homomorphism $\Pi:\Fr_n^\xi\to G/Z$
satisfying \eqref{eq:exp-exp}. Since $\Fr_n^\xi$ is simply
connected, we   can lift (in a unique way) $\Pi$ to a homomorphism 
$\wt\Pi: \Fr_n^\xi\to G$, see \cite{Pon}, Sect. 51, Theorem 80. 
A proof of  Theorem \ref{th0} is the same.
\hfill $\square$

\tt 

University of Graz,
\\
\phantom{.}
\hfill Department of Mathematics and Scientific computing;

High School of Modern Mathematics MIPT,
\\
\phantom{.}
\hfill 1 Klimentovskiy per., Moscow; 

Moscow State University, MechMath. Dept;

 University of Vienna, Faculty of Mathematics.
 
 \sm

e-mail:yurii.neretin(dog)univie.ac.at

URL: https://www.mat.univie.ac.at/$\sim$neretin/ 

\phantom{URL:} https://imsc.uni-graz.at/neretin/index.html

\end{document}